\documentstyle{amsppt}
\pagewidth{5.26in} \pageheight{7.86in} \NoRunningHeads
\magnification = 1200

\topmatter
\title {Comparison Theorems in pseudo-Hermitian Geometry and
applications}
\endtitle
\author Yuxin Dong and Wei Zhang \endauthor
\thanks {*Supported by NSFC grant No. 11271071, and LMNS, Fudan.}
\endthanks
\abstract{In this paper, we study the theory of geodesics with
respect to the Tanaka-Webster connection in a pseudo-Hermitian
manifold, aiming to generalize some comparison results in
Riemannian geometry to the case of pseudo-Hermitian geometry. Some
Hopf-Rinow type, Cartan-Hadamard type and Bonnet-Myers type
results are established.}
\endabstract
\endsubjclass
\endkeywords
\endtopmatter
\document
\heading{\bf Introduction}
\endheading
\vskip 0.3 true cm

A CR structure on an $(2m+1)$-dimensional manifold $M^{2m+1}$ is an $2m$%
-dimensional distribution $H(M)$ endowed with a formally
integrable complex structure $J$. The geometry of CR manifolds
goes back to Poincar\'e and received a great attention in the
works of Cartan, Tanaka, Chern-Moser, and others (cf. [Jo]). There
have been, over the last twenty or thirty years, many studies in
geometry and analysis on CR manifolds (cf. [DT], [BG], [CT], [VZ],
[CCY], [CKT]).

A pseudo-Hermitian manifold, which is an odd-dimensional analogue
of Hermitian manifolds, is a CR manifold $M$ endowed with a
pseudo-Hermitian structure $\theta $. The pseudo-Hermitian
structure $\theta $ determines uniquely a global nowhere zero
vector field $\xi $ and it, combining with the complex structure
$J$, induces a Riemannian metric $g_\theta $ on $M$ too. It turns
out that the Levi-Civita connection $\nabla^{\theta}$ of
$g_{\theta}$ is not convenient for investigating the
pseudo-Hermitian manifold, because it is not compatible with the
CR structure. From [Ta], [We], we know that there is a unique
canonical linear connection $\nabla $ (the Tanaka-Webster
connection), which is compatible with both the metric $g_\theta $
and the CR structure (see Proposition 2.1). This connection always
has nonvanishing torsion $T_\nabla (\cdot ,\cdot )$, whose partial
component $T_\nabla (\xi ,\cdot )$ is an important
pseudo-Hermitian invariant, called the pseudo-Hermitian torsion. A
Sasakian manifold, which is an odd dimensional analogue of
K\"ahler manifolds, is a pseudo-Hermitian manifold with vanishing
pseudo-Hermitian torsion. Besides its similarity with K\"ahler
geometry, interest in Sasakian manifolds has been from theoretical
physics with AdS/CFT correspondence, which provides a duality
between field theories and string theories.

In this paper, we investigate the theory of $\nabla $-geodesics on
the pseudo-Hermitian manifold aiming to generalize some comparison
results in Riemannian geometry, such as Hopf-Rinow type,
Cartan-Hadamard-type and Bonnet-Myers-type theorems, etc., to the
case of pseudo-Hermitian geometry. For this purpose, we shall
study the exponential maps, conjugate points and Jacobi fields
with respect to the Tanaka-Webster connection $\nabla$. The
organization of this paper is the following: In Section 1, we
recall some basic notions and properties of pseudo-Hermitian
manifolds. Section 2 is devoted to the exponential map $\exp
^\nabla $ and Hopf-Rinow type theorem. We find that the Gauss
lemma for $\nabla$-geodesics is no longer true due to the torsion
of $\nabla$. As a result, a short $\nabla$-geodesic is not
necessarily a length-minimizing curve. This causes some trouble
for establishing Hopf-Rinow type theorem. In order to study the
metric properties of $(M,\nabla)$, a natural distance $\delta$
between any two points is introduced by taking infimum of the
lengths of all broken $\nabla$-geodesics joining the two points.
In terms of the distance $\delta$, we establish a partial
Hopf-Rinow type result for pseudo-Hermitian manifolds, which
states that if $(M,\delta)$ is complete, then $(M,\nabla)$ is
complete. In Section 3, we investigate the Jacobi fields along a
$\nabla $-geodesic. As in the case for a Riemannian manifold, for
any two vector $v,w\in T_pM$, $(d\exp _{p}^{\nabla })_{tv}(tw)$ is
a Jacobi field along the geodesic $\gamma (t)=exp^{\nabla}_p(tv)$.
We compute the Taylor expansion of $\Vert (d\exp _{p}^{\nabla
})_{tv}(tw)\Vert ^{2}$ to show that the behavior of
$\nabla$-geodesics $exp^{\nabla}_p((v+sw)t)$ is affected by both
the curvature and the torsion of $(M,\nabla)$. Next, by
generalizing a result in [BD], the decomposition of a Jacobi field
$V(t)$ along any geodesic $\gamma (t)$ is given with respect to
$\gamma ^{\prime }(t)$ and its complementary space. Finally in
this section, we give explicitly the Jacobi fields along geodesics
in the Heisenberg group. In Section 4, we study Cartan-Hadamard
type result for pseudo-Hermitian manifolds. The main result in
this section asserts that if $M^{2m+1}$ is a complete Sasakian
manifold with non-positive horizontal curvature, then for any
$p\in M$, $exp^{\nabla}_p:T_pM\rightarrow N$ is a covering map.
Hence the covering space of $M$ is diffeomorphic to $R^{2m+1}$.
Finally, Section 5 is devoted to establishing a index comparison
theorem. As applications, we get a Bonnet-Myers type result
concerning conjugate points along geodesics in Sasakian manifolds
with either positive horizontal sectional curvature or positive
Ricci curvature. We should mention that various geodesics on
pseudo-Hermitian manifolds have been investigated by several
authors yet from somewhat different viewpoints. It is known that
the general theory of sub-Riemannian geodesics was established in
[St], which is a Hamiltonian description about geodesics in
cotangent bundles. Since a pseudo-Hermitian manifold may be
regarded as a special sub-Riemannian manifold, one may apply the
sub-Riemannian geodesic theory to pseudo-Hermitian manifolds.
Actually Barletta and Dragomir [BD] re-expressed the
sub-Riemannian geodesic equation via the Tanaka-Webster connection
$\nabla $ and studied the relationship between the sub-Riemannian
geodesics and the $\nabla$-geodesics on a pseudo-Hermitian
manifold. Besides, they also established some Cartan-Hadamard type
and Bonnet-Myers type results for conjugates along horizontal
$\nabla $-geodesics. In fact, the main results in [BD] involve
horizontal geodesics. By getting rid of the horizontal restriction
for geodesics, we are able to generalize some results in [BD] to
$\nabla$-geodesics with initial tangent vectors in any directions.
We shall find that although the connection $\nabla $ shares some
common notions and properties with general linear connections, it
displays some special features of a pseudo-Hermitian manifold too.
Furthermore, some interesting geometric properties of
pseudo-Hermitian manifolds are invisible from the Levi-Civita
connection $\nabla ^\theta $, but visible from the Tanaka-Webster
connection $\nabla $.

\heading{\bf 1. Preliminaries}
\endheading
\vskip 0.3 true cm

In this section, we  collect some facts and notations concerning
pseudohermitian structures on CR manifolds (cf. [DT], [BG] for
details).

\definition{Definition 1.1} Let $M^{2m+1}$ be a real $(2m+1)$-dimensional orientable $C^\infty $
manifold. A CR structure on $M$ is a complex rank $m$ subbundle
$H^{1,0}M$ of $TM\otimes C$ satisfying
\newline
(i) $H^{1,0}M\cap H^{0,1}M=\{0\}$
($H^{0,1}M=\overline{H^{1,0}M}$);
\newline
(ii) $[\Gamma (H^{1,0}M),\Gamma (H^{1,0}M)]\subseteq \Gamma
(H^{1,0}M)$.
\newline The pair $(M,H^{1,0}M)$ is called a CR manifold.
\enddefinition

The complex subbundle $H^{1,0}M$ corresponds to a real subbundle
of $TM$:
$$
H(M)=Re\{H^{1,0}M\oplus H^{0,1}M\}\tag{1.1}
$$
which is endowed with a natural complex structure $J$ as follows
$$
J(V+\overline{V})=i(V-\overline{V})\tag{1.2}
$$
for any $V\in H^{1,0}M$. Equivalently, the CR structure may be
described by the pair $(H(M),J)$. A $C^\infty $ map
$f:(M,H(M),J)\rightarrow (N,H(N), \widetilde{J})$ between two CR
manifolds is called a CR map if $df(H(M))\subset H(N)$ and
$(df\circ J)|_{H(M)}=(\widetilde{J}\circ df)|_{H(M)}$.

Let $E$ be the conormal bundle of $H(M)$ in $T^{*}M$, whose fiber
at each point $x\in M$ is given by
$$E_x=\{\omega \in T_x^{*}M:\ker
\omega \supseteq H_x(M)\}. \tag{1.4}
$$
Since both $TM$ and $H(M)$ are orientable, $E\simeq T(M)/H(M)$ is
an orientable line bundle. It is known that any orientable real
line bundle over a connected manifold is trivial. Therefore $E$
admits globally defined nowhere vanishing sections.

\definition{Definition 1.2} A globally defined nowhere vanishing section $\theta \in \Gamma
(E)$ is called a pseudo-Hermitian structure on $M$. The Levi-form
$L_\theta $ associated with a pseudo-Hermitian structure $\theta$
is defined by
$$
L_\theta (X,Y)=d\theta (X,JY)\tag{1.4}
$$
for any $X,Y\in H(M)$. If $L_\theta $ is positive definite for
some $\theta$, then $(M,H(M),J)$ is said to be strictly
pseudoconvex.
\enddefinition

Henceforth we assume that $(M,H(M),J)$ is a strictly pseudoconvex
CR manifold endowed with a pseudo-Hermitian structure $\theta $
such that $L_\theta $ is positive definite. The quadruple
$(M,H(M),J,\theta )$ is called a pseudo-Hermitian manifold, which
is sometimes denoted simply by $(M,\theta )$.

For a pseudo-Hermitian manifold $(M,H(M),J,\theta )$, there exists
a unique globally defined nowhere zero vector field $\xi $ such
that (cf.[DT])
$$
\theta (\xi )=1,\quad d\theta (\xi ,\cdot )=0.\tag{1.5}
$$
This vector field $\xi $ is referred to as the Reeb vector field,
which is transverse to $H(M)$. Consequently we have the following
decomposition
$$
TM=L_{\xi}\oplus H(M) ,\tag{1.6}
$$
where $L_{\xi}$ is the trivial line bundle generated by $\xi$. The
subbundles $L_{\xi}$ and $H(M)$ will be called the vertical and
horizontal distributions respectively. Correspondingly, a vector
$V\in TM$ is called vertical (resp. horizontal) if $V\in L_{\xi}$
(resp. $H(M)$). For convenience, we extend $J$ to a $(1,1)$-tensor
field on $M$ by requiring that
$$
J\xi =0.\tag{1.7}
$$

Let $\pi _H:TM\rightarrow H(M)$ be the natural projection
morphism. Set $G_\theta =\pi _H^{*}L_\theta $, that is,
$$
G_\theta (X,Y)=L_\theta (\pi _HX,\pi _HY)\tag{1.8}
$$
for any $X$, $Y\in TM$. Then one may introduce a Riemannian
metric, called the Webster metric, as follows
$$
g_\theta =\theta \otimes \theta +G_\theta\tag{1.9}
$$
which is sometimes denoted by $\langle \cdot ,\cdot \rangle $ for
simplicity. Clearly
$$
\theta (X)=\langle \xi ,X\rangle ,\ d\theta (X,Y)=\langle
JX,Y\rangle .\tag{1.10}
$$
In terms of (1.9) and (1.10), we find that (1.6) is actually an
orthogonal decomposition. In addition, $\theta \wedge (d\theta
)^m$ is, up to a constant, the volume form of $(M,g_\theta )$.

On a pseudo-Hermitian manifold, we have the following canonical
linear connection which preserves both the CR and the metric
structures.

\proclaim{Proposition 1.1 ([Ta], [We])} Let $(M,H(M),J,\theta )$
be a pseudo-Hermitian manifold. Then there exists a unique linear
connection $\nabla$ such that
\newline(i) $\nabla _X\Gamma
(H(M))\subset \Gamma (H(M))$ for any $X\in \Gamma (TM)$;
\newline (ii) $\nabla g_\theta =0$, $\nabla J=0$ (hence $\nabla
\xi=\nabla \theta =0$);
\newline (iii) The torsion $T_\nabla $ of
$\nabla $ is pure, that is, for any $X,Y\in H(M)$, $T_\nabla
(X,Y)=d\theta (X,Y)\xi$ and $T_\nabla (\xi,JX)+JT_\nabla
(\xi,X)=0$.
\endproclaim

The connection $\nabla$ in Proposition 1.1 is called the
Tanaka-Webster connection. Note that the torsion of the
Tanaka-Webster connection is always non-zero. The pseudo-Hermitian
torsion, denoted by $\tau$, is the $TM$-valued $1$-form defined by
$\tau (X)=T_\nabla(\xi,X)$. The anti-symmetry of $T_\nabla$
implies that
$$
\tau(\xi)=0.\tag{1.11}
$$
Using (iii) of Proposition 1.1 and the definition
of $\tau$, the total torsion of the Tanaka-Webster connection may
be expressed as
$$
T_\nabla (X,Y)=2(\theta \wedge \tau )(X,Y)+2d\theta (X,Y)\xi
\tag{1.12}
$$
for any $X,Y\in TM$. Set
$$
A(X,Y)=g_\theta (\tau X,Y)  \tag{1.13}
$$
for any $X,Y\in TM$. Then the properties of $\nabla$ in
Proposition 1.1 also imply that $\tau(H^{1,0}(M))\subset
H^{0,1}(M)$ and $A$ is a trace-free symmetric tensor field.

\proclaim{Lemma 1.2 (cf. [DT])} The Levi-Civita connection $\nabla
^\theta $ of $(M,g_\theta )$ is related to the Tanaka-Webster
connection by
$$
\nabla ^\theta =\nabla -(d\theta +A)\xi+\tau \otimes \theta
+2\theta \odot J
$$
where $(\theta \odot J)(X,Y)=\frac 12(\theta (X)JY+\theta (Y)JX)$
for any $ X,Y\in TM$.
\endproclaim

For a pseudo-Hermitian manifold $(M,\theta )$, the curvature
tensor $R$ with respect to its Tanaka-Webster connection is
defined by
$$
R(X,Y)Z=\nabla _X\nabla _YZ-\nabla _Y\nabla _XZ-\nabla _{[X,Y]}Z
$$
for any $X,Y,Z\in \Gamma (TM)$. Clearly $R$ satisfies
$$
\langle R(X,Y)Z,W\rangle =-\langle R(Y,X)Z,W\rangle =-\langle
R(X,Y)W,Z\rangle\tag{1.14}
$$
where the second equality is because of $\nabla g_\theta =0$.
However, the symmetric property $\langle R(X,Y)Z,W\rangle =\langle
R(Z,W)X,Y\rangle $ is no longer true for a general
pseudo-Hermitian manifold due to the failure of the first Bianchi
identity.

For a horizontal $2$-plane $P=span_R\{X,Y\}\subset H(M)$, the
horizontal sectional curvature of $P$ is defined by
$$
K^H(P)=\frac{\langle R(X\wedge Y),X\wedge Y\rangle}{\langle
X\wedge Y,X,\wedge Y\rangle}.\tag{1.15}
$$
We define the Ricci tensor of $\nabla $ by
$$
Ric(Y,Z)=trace\{X\mapsto R(X,Z)Y\}\tag{1.16}
$$
for any $Y,Z\in TM$.

\definition{Definition 1.3} A pseudo-Hermitian manifold $(M,H(M),J,\theta)$ is called a
Sasakian manifold if its pseudo-Hermitian torsion $\tau$ is zero.
\enddefinition

From [DT], we know that if $(M,\theta )$ is a Sasakian manifold,
then
$$
\langle R(X,Y)Z,W\rangle =\langle R(Z,W)X,Y\rangle .\tag{1.17}
$$
Consequently, if one of $X,Y,Z$ and $W$ is vertical, then
$$
\langle R(X,Y)Z,W\rangle =0.\tag{1.18}
$$

We denote the length of a continuous piecewise smooth curve $
c:[a,b]\rightarrow (M,g_\theta )$ by $L[c]$.

\proclaim{Lemma 1.3} (see also \S 7 in [BD]) Let $\alpha
:[a,b]\times (-\varepsilon ,\varepsilon )\rightarrow (M,\theta )$
be a smooth map. Set $T=d\alpha (\frac \partial {\partial t})$,
$V=d\alpha (\frac
\partial {\partial s})$. Write $\alpha (t,s)$ as $\alpha _s(t)$
($a\leq t\leq b$, $-\varepsilon <s<\varepsilon $). Assume that
$\alpha _0(t)$ is parameterized by arc length. Then
$$
\frac d{ds}L(\alpha _s)|_{s=0}=\langle V,T\rangle
|_a^b-\int_a^b\{\langle V,\nabla _TT\rangle -\langle T_\nabla
(V,T),T\rangle \}dt
$$
\endproclaim
\demo{Proof} The usual first variation formula of arc length gives
(cf. [CE])
$$
\frac d{ds}L(\alpha _s)|_{s=0}=\langle V,T\rangle
|_a^b-\int_a^b\langle V,\nabla _T^\theta T\rangle dt.\tag{1.19}
$$
The lemma follows immediately from Lemma 1.2, (1.12) and
(1.19).\qed
\enddemo

\heading{\bf 2. Exponential map and Hopf-Rinow Type Results}
\endheading
\vskip 0.3 true cm

Let $(M,H(M),J,\theta )$ be a pseudo-Hermitian manifold with the
Tanaka-Webster connection $\nabla $. A $C^1$ curve $\gamma
:[0,l]\rightarrow M$ is called a $\nabla $-geodesic if $\nabla
_{\gamma ^{\prime }}\gamma ^{\prime }=0$ on $[0,l]$. Since the
linear connection $\nabla $ is of class $ C^\infty $, a
$C^1$-geodesic with respect to $\nabla $ is automatically of class
$C^\infty $ (cf. [KN1]). A smooth curve $\gamma :[0,l]\rightarrow
M$ is referred to as a slant curve if the angle between $\gamma
^{\prime }(t)$ and $\xi _{\gamma (t)}$ is constant along $\gamma
$. In particular, if $ \gamma ^{\prime }(t)\ $is perpendicular
(resp. parallel) to $\xi _{\gamma (t)}$ for each $t$, then $\gamma
$ is called a horizontal (resp. vertical) curve. Since $\nabla
g_\theta =0$ and $\nabla \xi =0$, it is clear that any $\nabla
$-geodesic $\gamma (t)$ must be a slant curve. In particular, if
the initial tangent vector of the $\nabla $-geodesic is horizontal
(resp. vertical), then $\gamma $ should be horizontal (resp.
vertical).

Given a point $p\in M$ and a vector $v\in T_pM$, the ODE theory
implies that there exists a unique $\nabla $-geodesic $\gamma
_v(t)$ satisfying $\gamma _v(0)=p$ and $\gamma _v^{\prime }(0)=v$.
Since a parameterization which makes $\gamma $ into a geodesic, if
any, is determined up to an affine transformation of $t$, the
parameter $t$ is called an affine parameter. As usual, the
exponential map $\exp _p^\nabla :T_pM\rightarrow M$ is defined by
$$
\exp _p^\nabla (v)=\gamma _v(1)
$$
for all $v\in T_pM$ such that $1$ is in the domain of $\gamma _v$.
Since $\nabla $ preserves the metric $g_\theta $, we find that
$|\gamma _v^{\prime }(t)|=|v|$ for each $t$ and thus $L(\gamma
_v;[0,1])=|v|$, where $L(\cdot )$ denotes the length of the curve.
The linear connection $\nabla $ of $M$ is said to be complete if
every $\nabla $-geodesic can be extended to a geodesic $\gamma
(t)$ defined for $-\infty <t<\infty $, where $t$ is an affine
parameter. Hence, if $(M,\nabla )$ is complete, then $\exp ^\nabla
$ is defined on all of $TM$ and vice versa.

From Proposition 8.2 in [KN1], we know that there is a
neighborhood $D_p$ of each point $p$ (more precisely, the zero
vector at $p$) in $T_pM$ which is mapped diffeomorphically onto a
neighborhood $U_p$ of $p$ in $M$ by the exponential map. Choosing
a linear frame $u=\{X_1,...,X_n\}$ at $p$, the diffeomorphism
$\exp _p^\nabla :D_p\rightarrow U_p$ defines a local coordinate
system in $U_p$ in a natural manner. This local coordinate system
is called a normal coordinate system at $p$. The following result
holds for any linear connection on a manifold.

\proclaim{Proposition 2.1} (cf. [KN1]) Let $x^1,...,x^n$ be a
normal coordinate system with origin $p$. Let $U(p;\rho )$ be the
neighborhood of $p$ defined by $\sum (x^i)^2<\rho ^2$. Then there
is a positive number $a$ such that if $0<\rho <a$, then
\newline (1) $U(p;\rho )$ is convex in the sense that any two points of
$U(p;\rho )$ can be joined by a $\nabla $-geodesic which lies in
$U(p;\rho )$;
\newline (2) Each point of $U(p;\rho )$ has a normal coordinate
neighborhood containing $U(p;\rho )$.
\endproclaim

For each $v\in T_pM$, the tangent space $T_v(T_pM)$ can be
identified with $T_pM$ in a natural way. Due to the Webster
metric, each tangent space $T_pM$ comes equipped with an inner
product. Therefore $T_v(T_pM)$ inherits an inner product from
$(g_\theta )_p(\cdot ,\cdot )$.

\proclaim{Proposition 2.2} Let $c(s)$ be a curve in $T_pM$ such
that every point of $c$ is at the same distance from the origin of
$T_pM$. Let $\rho _s(t):[0,1]\rightarrow T_pM$ be the ray from $0$
to $c(s)$ in $T_pM$. Set $\alpha (t,s)=\exp (\rho _s(t))$ and
$\gamma (t)=\alpha (t,0)$. Assume that $\exp ^\nabla $ is defined
through $\rho _s$. Then
$$
\langle V(t),\gamma ^{\prime }(t)\rangle =\int_0^t[2\langle
J\gamma ^{\prime },V\rangle \theta (\gamma ^{\prime })+\theta
(\gamma ^{\prime })\langle \tau (\gamma ^{\prime }),V\rangle
-\theta (V)\langle \tau (\gamma ^{\prime }),\gamma ^{\prime
}\rangle ]\tag{2.1}
$$
where $V(t)=$ $d\alpha (\frac \partial {\partial s})_{(t,0)}$ is
the variation vector field of $\alpha $ along the $\nabla
$-geodesic $\gamma $.
\endproclaim
\demo{Proof} From the definition of $\exp ^{\nabla}$ and $\alpha
$, we know that the lengths of the curve $t\mapsto \alpha (t,s)$
are independent of $s$. Note that $V(0)=0$. Then Lemma 1.3 and
(1.12) yield
$$
\aligned 0&=\frac{dL[\exp (\rho
_s)|_{[0,t]}]}{ds}|_{s=0}\\
&=\langle V(t),\gamma ^{\prime }(t)\rangle -\int_0^t[2\theta
(\gamma ^{\prime })\langle J\gamma ^{\prime },V\rangle +\theta
(\gamma ^{\prime })\langle \tau (V),\gamma ^{\prime }\rangle
-\theta (V)\langle \tau (\gamma ^{\prime }),\gamma ^{\prime
}\rangle ]dt.\endaligned
$$
Consequently
$$
\langle V(t),\gamma ^{\prime }(t)\rangle =\int_0^t[2\langle
J\gamma ^{\prime },V\rangle \theta (\gamma ^{\prime })+\theta
(\gamma ^{\prime })\langle \tau (\gamma ^{\prime }),V\rangle
-\theta (V)\langle \tau (\gamma ^{\prime }),\gamma ^{\prime
}\rangle ].
$$
\qed\enddemo

\remark{Remark 2.1} For $v\in T_pM$, we assume that $w\in
T_v(T_pM)$ is perpendicular to $v$ when $w$ is also regarded as a
vector in $T_pM$. Clearly there exists a curve $c(s)$ in $T_pM$
such that $c(0)=v$, $c^{\prime }(0)=w$ and such that every point
of $c$ is at the same distance from the origin of $T_pM$. In this
circumstance, (2.1) becomes
$$
\aligned
&\langle (d\exp ^\nabla )_{tv}(tw),(d\exp ^\nabla )_{tv}(v)\rangle \\
&=\int_0^t[2\langle J\gamma ^{\prime },V\rangle \theta (\gamma
^{\prime })+\theta (\gamma ^{\prime })\langle \tau (\gamma
^{\prime }),V\rangle -\theta (V)\langle \tau (\gamma ^{\prime
}),\gamma ^{\prime }\rangle ].\endaligned\tag{2.2}
$$
In particular, if $M$ is Sasakian, then
$$
\langle (d\exp ^\nabla )_{tv}(tw),(d\exp ^\nabla )_{tv}(v)\rangle
=2\int_0^t\langle J\gamma ^{\prime },V\rangle \theta (\gamma
^{\prime }).\tag{2.3}
$$
\endremark

From either (2.2) or (2.3), we see that the Gauss lemma is no
longer true for $\exp ^\nabla $ due to the torsion. However, it
still holds for some special geodesics.

\proclaim{Corollary 2.3} Let $\rho (t)=tv$ ($t\in [0,1]$) be a ray
through the origin in $T_pM$ and let $w\in T_v(T_pM)$ be a vector
perpendicular to $v$. Assume that $\exp ^\nabla $ is defined along
$\rho $. If either $v$ is vertical or $M$ is Sasakian and $v$ is
horizontal, then we have $$\langle (d\exp ^\nabla
)_{tv}(tw),(d\exp ^\nabla )_{tv}(v)\rangle =0$$ for each $t\in
[0,1]$.
\endproclaim
\demo{Proof} First, assume that $v$ is vertical. Then $\gamma
^{\prime }(t)$ is vertical for each $t$. In terms of (1.7) and
(1.11), we deduce from (2.2) that $\langle (d\exp ^\nabla
)_{tv}(tw),(d\exp ^\nabla )_{tv}(v)\rangle =0$.

Next, assume that $v\in H_p(M)$. Then $\gamma (t)=\exp ^\nabla
\rho (t)$ is a horizontal geodesic. Furthermore, if $M$ is
Sasakian, we get from immediately (2.3) the required result.\qed
\enddemo

In order to investigate the metric properties of the $\nabla
$-geodesics, we shall introduce a distance function determined by
the connection $\nabla $ and $g_\theta $. A continuous curve
$c:[a,b]\rightarrow M$ is called a broken $\nabla $-geodesic if
there exist $a=a_1<a_2<\cdot \cdot \cdot <a_n=b$ such that
$c:[a_i,a_{i+1}]\rightarrow M$ is a $\nabla $-geodesic for $
i=1,...,n-1$. Clearly a broken $\nabla $-geodesic is piecewise
$C^\infty $ -differentiable. For any $p,q\in M$, let $\Gamma
(p,q)\ $denote all broken $\nabla $-geodesics joining $p$ and $q$.
We define the distance between $p$ and $q$, $\delta (p,q)$, by
$$\delta (p,q)=\inf_{\gamma \in \Gamma (p,q)}L(\gamma ).\tag{2.4}
$$
For any continuous curve connecting $p$ and $q$, it can be covered
by a finite number of normal coordinate neighborhoods. Thus there
always exist broken $\nabla $-geodesics joining the two points, so
the distance $\delta $ is finite. It is easy to verify that
$(M,\delta )$ is a metric space. Let $d$ denote the usual
Riemannian distance function determined by $g_\theta $. Clearly
$d(p,q)\leq \delta (p,q)$ for any $p,q\in M$. Hence if $(M,d)$ is
complete, then $(M,\delta )$ is complete too.

\proclaim{Proposition 2.4} The distance function $\delta $ defines
the same topology as the manifold topology of $M$.
\endproclaim

\demo{Proof} It is known that the Riemannian distance function $d$
defines the same topology as the original topology of $M$ (cf.
[KN1], [Ch]). Therefore we only need to verify that $d$ and
$\delta $ define the same metric space topology.

Suppose first that $U$ is an open subset of the metric space
topology defined by $d$. Since $d(p,q)\leq \delta (p,q)$ for any
$p,q\in M$, we find that $U$ must be an open subset of the metric
space topology defined by $\delta $.

Assume now that $W$ is an open subset defined by $\delta $, that
is, for any point $p\in W$, there exits $\varepsilon >0$ such that
$B_\delta (p;\varepsilon )\subset W$, where $B_\delta
(p;\varepsilon )=\{q:\delta (p,q)<\varepsilon \}$. Set
$$
D_p(\varepsilon )=\{v\in T_pM:||v||<\varepsilon \},\quad B^\nabla
(p;\varepsilon )=\exp _p^\nabla (D_p(\varepsilon )).
$$
Let $\exp ^{\nabla ^\theta }$ denote the Riemannian exponential
map. For sufficiently small $\varepsilon $, we know that $\exp
_p^{\nabla ^\theta }(D_p(\varepsilon ))=B_d(p;\varepsilon )$ (cf.
[KN1], [Ch]). Moreover, both $\exp _p^\nabla :D_p(\varepsilon
)\rightarrow B^\nabla (p;\varepsilon )$ and $\exp _p^{\nabla
^\theta }:D_p(\varepsilon )\rightarrow B_d(p;\varepsilon )$ are
diffeomorphisms. Consequently $B^\nabla (p;\varepsilon )$ is an
open subset in $(M,d)$. Note that if $\gamma $ is a
$\nabla$-geodesic joining $p$ to a point $q$ in the normal
coordinate neighborhood $B^\nabla (p;\varepsilon )$, then $\delta
(p,q)\leq L(\gamma )<\varepsilon $. This implies that $B^\nabla
(p;\varepsilon )\subset B_\delta (p;\varepsilon )$, and thus
$B^\nabla (p;\varepsilon )\subset W$. We then conclude that $W$ is
also an open subset defined by $d$. \qed
\enddemo

\proclaim{Theorem 2.5} Let $(M,H(M),J,\theta )$ be a
pseudo-Hermitian manifold. If $(M,\delta )$ is complete, then
$(M,\nabla )$ is complete, that is, $\exp ^\nabla $ is defined on
all of $TM$.
\endproclaim
\demo{Proof} Given $p\in M$ and $0\neq v\in T_pM$, we assume that
$\gamma $ is a $\nabla $-geodesic with $\gamma (0)=p$ and $\gamma
^{\prime }(0)=v$. Suppose $[0,t_0)$ is the largest open interval
for which such a $\gamma $ exists. Note that the parameter $t$
must be proportional to arc length. Hence, if $t_0$ is finite and
$t_i\nearrow t_0$, then $\delta (\gamma (t_i),\gamma (t_j))\leq
L(\gamma |_{[t_i,t_j]})=c|t_j-t_i|\rightarrow 0$ as
$i,j\rightarrow \infty $ ($i<j$) for some positive constant $c$.
Consequently $\{\gamma (t_i)\}$ is a Cauchy sequence with some
limit $q$ in $(M,\delta )$. Define $\gamma (t_0)=q$. Let $U(q;\rho
)$ be a normal coordinate system as in Proposition 2.1. For
sufficiently large $i$, $\gamma (t_i)\in U(q;\rho )$. Let $\sigma
:[0,r_0)\rightarrow M$ be a $\nabla $-geodesic with $\sigma
(0)=\gamma (t_i)$ and $\sigma ^{\prime }(0)=\gamma ^{\prime
}(t_i)$ and let $[0,r_0)$ be the largest open interval for which
$\sigma (t)$ exists. According to Proposition 2.1, $\gamma (t_i)$
has a normal coordinate neighborhood containing $U(q;\rho )$. Thus
$r_0>t_0-t_i$ and $\gamma (t_0)\in \sigma $. Therefore $ \gamma
\cup \sigma $ is a smooth $\nabla $-geodesic, and $\gamma $
extends past $t_0$, which is a contradiction.\qed
\enddemo

The Hopf-Rinow theorem in Riemannian geometry tells us that the
completeness of $(M,\nabla ^\theta )$ is equivalent to the
completeness of $(M,d)$. Since the completeness of $(M,d)$ yields
the completeness of $(M,\delta )$, we have

\proclaim{Corollary 2.6} Let $(M,H(M),J,\theta )$ be a
pseudo-Hermitian manifold. If $(M,\nabla ^\theta )$ is complete,
then $(M,\nabla )$ is complete.
\endproclaim

\heading{\bf 3. Jacobi fields on pseudo-hermitian manifolds}
\endheading
\vskip 0.3 true cm

From now on, we always assume that $(M,\nabla )$ is a complete
pseudo-hermitian manifold of dimension $2m+1$. Let us consider a
$1$-parameter family of $\nabla $-geodesics given by a map $\alpha
(t,s):[a,b]\times (-\varepsilon ,\varepsilon )\rightarrow M$ such
that for each fixed $s$, $\alpha (t,s)$ is a $\nabla $-geodesic.
Set $T=d\alpha (\frac \partial {\partial t})$ and $V=d\alpha
(\frac \partial {\partial s})$. Since $[\partial /\partial
t,\partial /\partial s]=0$, we have
$$
\nabla _TV-\nabla _VT=T_\nabla (T,V).\tag{3.1}
$$
Therefore
$$
\nabla _T\nabla _TV=\nabla _T\nabla _VT+\nabla _T(T_\nabla
(T,V)).\tag{3.2}
$$
By the definition of curvature tensor and the geodesic equation
$\nabla _TT=0$, it follows from (3.2) that (see also Theorem 1.4
in [Ch])
$$
\nabla _T\nabla _TV=R(T,V)T+\nabla _T(T_\nabla (T,V)).\tag{3.3}
$$
The equation above is called the Jacobi equation. A vector field
$V$ satisfying the equation (3.3) is called a Jacobi field  along
the geodesic $\gamma$. For example, it is easy to verify that
$(a+bt)\gamma ^{\prime }$ is a Jacobi field for any $a,b\in R$.
However, a general solution of the Jacobi equation can not be
given so explicitly and should depend on both the curvature and
the torsion tensors. Since (3.3) is an ODE system of second order,
a Jacobi field $V$ is determined uniquely by $V(0)$ and $V^{\prime
}(0)$. Let $J_\gamma $ denote the real linear space of all Jacobi
fields along $\gamma $ in $(M,\nabla )$. Then $\dim _RJ_\gamma
=4m+2$. We have shown that the variation field of a $1$-parameter
family of geodesics is a Jacobi field. Conversely, if $V$ is a
Jacobi field, then $V$ comes from a variation of geodesics too. In
fact, let $c(s)$ be a curve such that $c^{\prime }(0)=V(0)$, and
let $\gamma ^{\prime }(0)$, $V(0)$ and $T_\nabla (\gamma ^{\prime
}(0),V(0))$ be extended respectively to parallel fields $\gamma
^{\prime }(0)_s$, $V^{\prime }(0)_s$ and $T_\nabla (\gamma
^{\prime }(0),V(0))_s$ along $c(s)$. Then the variation field of
$$
\alpha (t,s):=\exp _{c(s)}^\nabla \{t[\gamma ^{\prime
}(0)_s+sV^{\prime }(0)_s-sT_\nabla (\gamma ^{\prime
}(0),V(0))_s]\}
$$
is a Jacobi field $\widetilde{V}(t)=d\alpha (\frac \partial
{\partial s})|_{s=0}$. Clearly $\widetilde{V}(0)=V(0)$. Using
(3.1), we compute
$$
\aligned \widetilde{V}^{\prime }(0) &=\{\nabla _{\frac \partial
{\partial t}}d\alpha
(\frac \partial {\partial s})|_{s=0}\}|_{t=0} \\
\ &=\nabla _{\frac \partial {\partial s}}\{[\gamma ^{\prime
}(0)_s+sV^{\prime }(0)_s-sT_\nabla (\gamma ^{\prime
}(0),V(0))_s]\}|_{s=0}+T_\nabla (\gamma ^{\prime }(0),V(0)) \\
\ &=V^{\prime }(0).\endaligned
$$
Since $\widetilde{V}(t)$ and $V(t)$ have the same initial
conditions, we conclude that $\widetilde{V}(t)=V(t)$. Consequently
a Jacobi field $V$ with $V(0)=0$ and $V^{\prime }(0)=w$ may be
given by
$$
V(t)=\frac \partial {\partial s}\exp _p^\nabla \left(
t(v+sw)\right) |_{s=0}=(d\exp _p^\nabla )_{tv}(tw).
$$
If we put $v=w$, then $V(t)=t\gamma ^{\prime }(t)$, which implies
$$
|(d\exp _p^\nabla )_v(v)|=|v|
$$
or in other words, that $\exp _p^\nabla :T_pM\rightarrow M$ is an
isometry in the radial direction.

We would like to obtain information on $\Vert (d\exp _{p}^{\nabla
})_{tv}(tw)\Vert ^{2}$ by calculating its Taylor expansion, where
$v,w\in T_{p}M$ are two unit vectors with $\langle v,w\rangle =0$.
Let $\gamma (t)=\exp _{p}^{\nabla }(tv)$ and $V(t)=(d\exp
_{p}^{\nabla })_{tv}(tw)$. From the above discussion, we know that
$V$ is a Jacobi field along $\gamma $ satisfying $V(0)=0$,
$V^{\prime }(0)=w$. Then
$$
\aligned &\langle V,V\rangle |_{t=0}=0, \quad
\langle V,V\rangle ^{\prime }=2\langle V,V^{\prime }\rangle |_{t=0}=0, \\
&\langle V,V\rangle ^{\prime \prime }|_{t=0}=2\langle V^{\prime
},V^{\prime }\rangle |_{t=0}+2\langle V^{\prime \prime },V\rangle
|_{t=0}=2.
\endaligned
$$
Note that $V^{\prime \prime }|_{t=0}=R(v,V(0))v+(T_{\nabla
}(\gamma ^{\prime },V))^{\prime }|_{t=0}=T_{\nabla }(v,w)$, so
$$
\aligned \langle V,V\rangle ^{\prime \prime \prime
}|_{t=0}=&6\langle V^{\prime \prime },V^{\prime }\rangle
|_{t=0}+2\langle V^{\prime \prime \prime },V\rangle
|_{t=0} \\
&=6\langle T_{\nabla }(v,w),w\rangle.\endaligned
$$
Also
$$
\aligned V^{\prime \prime \prime }|_{t=0}&=(\nabla _{\gamma
^{\prime }}R)(\gamma ^{\prime },V)\gamma ^{\prime
}|_{t=0}+R(\gamma ^{\prime },V^{\prime })\gamma ^{\prime
}|_{t=0}+[T_{\nabla }(\gamma ^{\prime },V)]^{\prime \prime
}|_{t=0}\\
&=R(v,w)v+[T_{\nabla }(\gamma ^{\prime },V)]^{\prime \prime
}|_{t=0}.\endaligned
$$
Then
$$
\aligned \langle V,V\rangle ^{\prime \prime \prime \prime
}&=8\langle V^{\prime \prime \prime },V^{\prime }\rangle
|_{t=0}+6\langle V^{\prime \prime },V^{\prime \prime }\rangle
|_{t=0}+2\langle V^{\prime \prime \prime \prime },V\rangle
|_{t=0} \\
&=8\langle R(v,w)v+[T_{\nabla }(\gamma ^{\prime
},V)]_{t=0}^{\prime \prime },w\rangle +6\langle T_{\nabla
}(v,w),T_{\nabla }(v,w)\rangle.\endaligned
$$
Therefore
$$
\aligned &\Vert d\exp ^{\nabla }(tw)\Vert^{2}=t^{2}+\langle
T_{\nabla }(v,w),w\rangle t^{3}\\
&+\frac{1}{4!}\{8\langle R(v,w)v+[T_{\nabla }(\gamma ^{\prime
},V)]_{t=0}^{\prime \prime },w\rangle +6\langle T_{\nabla
}(v,w),T_{\nabla }(v,w)\rangle
\}t^{4}+O(t^{5}).\endaligned\tag{3.4}
$$
In particular, if $(M,\theta )$ is Sasakian, then (3.4) and (2.12)
imply that
$$
\aligned \Vert d\exp ^{\nabla }(tw)\Vert^{2}&=t^{2}+2\langle
Jv,w\rangle \langle \xi
,w\rangle t^{3}-\frac{1}{3}\langle R(w,v)v,w\rangle t^{4} \\
&+\langle Jv,w\rangle ^{2}t^{4}+O(t^{5}).\endaligned\tag{3.5}
$$
Set $\rho _{s}(t)=(v+sw)t$. We discover from either (3.4) or (3.5)
that the behavior of geodesics $\exp _{p}(\rho _{s})$ is affected
by both the curvature and the torsion of $M$. Let us check some
special cases of (3.5) on a Sasakian manifold. For example, if $v$
is vertical, then $\exp ^{\nabla }$ nearly preservers the "width"
between the ray $\rho _{0}$ and $\rho _{s}$ with the error term
$O(t^{5})$. Next, if $v$ $w$ are both horizontal and $\langle
R(w,v)v,w\rangle $ is negative, then the geodesics locally diverge
when compared to the rays $\rho _{s}$. Finally, if $v, w$ are both
horizontal and satisfy the additional condition $\langle
Jv,w\rangle =0$, then the expansion (3.5) is almost same as that
for Riemannian exponential map (cf. [CE]). In this circumstance,
if $\langle R(w,v)v,w\rangle $ is positive, then the corresponding
geodesics locally converge by comparison with the rays $\rho
_{s}$.

Let $V$ be a Jacobi field along a $\nabla $-geodesic $\gamma $ in
a pseudo-Hermitian manifold $M$. Using (2.12) and the properties
that $\nabla d\theta =0$, $\nabla \theta =0$, we compute
$$
\aligned &\nabla _{\gamma ^{\prime }}T_\nabla (\gamma ^{\prime
},V)=\nabla _{\gamma ^{\prime }}\{2d\theta (\gamma ^{\prime
},V)\xi +2(\theta \wedge \tau
)(\gamma ^{\prime },V)\} \\
&=2d\theta (\gamma ^{\prime },\nabla _{\gamma ^{\prime }}V)\xi
+\theta (\gamma ^{\prime })\nabla _{\gamma ^{\prime }}\tau
(V)-\theta (\nabla _{\gamma ^{\prime }}V)\tau (\gamma ^{\prime
})-\theta (V)\nabla _{\gamma ^{\prime }}\tau (\gamma ^{\prime
}).\endaligned\tag{3.6}
$$
Consequently $V$ satisfies
$$
\aligned \nabla _{\gamma ^{\prime }}\nabla _{\gamma ^{\prime
}}V&=R(\gamma ^{\prime },V)\gamma ^{\prime }+2d\theta (\gamma
^{\prime },\nabla _{\gamma ^{\prime
}}V)\xi +\theta (\gamma ^{\prime })\nabla _{\gamma ^{\prime }}\tau (V) \\
&-\theta (\nabla _{\gamma ^{\prime }}V)\tau (\gamma ^{\prime
})-\theta (V)\nabla _{\gamma ^{\prime }}\tau (\gamma ^{\prime
}).\endaligned\tag{3.7}
$$
Note that $\langle \nabla _{\gamma ^{\prime }}\nabla _{\gamma
^{\prime }}V,\gamma ^{\prime }\rangle \ $is not necessarily
vanishing. Therefore, unlike the Riemannian case, the tangential
component $\langle V,\gamma ^{\prime }\rangle $ is not linear in
general. Although $(a+bt)\gamma ^{\prime }$ is a Jacobi field for
any $a$, $b\in R$, the tangential component of a Jacobi field may
contain some nonlinear part. We shall investigate this nonlinear
tangential part of a Jacobi field on a Sasakian manifold.

Suppose now that $M$ is a Sasakian manifold. Then (3.7) is
simplified to
$$
\nabla _{\gamma ^{\prime }}\nabla _{\gamma ^{\prime }}V=R(\gamma
^{\prime },V)\gamma ^{\prime }+2\langle J\gamma ^{\prime },\nabla
_{\gamma ^{\prime }}V\rangle \xi .\tag{3.8}
$$
First let us consider a vertical Jacobi field $V$. Writing
$V(t)=f(t)\xi _{\gamma (t)}$ and substituting it into (3.8), we
get
$$
f^{\prime \prime }=0,
$$
that is, $f(t)=a+bt$ for some $a,b\in R$. Hence we find that any
vertical Jacobi field along $\gamma $ is of the form $(a+bt)\xi
_{\gamma (t)}$ for some $a,b\in R$.

\proclaim{Lemma 3.1} Let $M$ be a Sasakian manifold and let
$\gamma $ be a $\nabla $-geodesic in $M$. Then any Jacobi field
$V$ along $\gamma $ satisfies
$$
\frac d{dt}\langle V,\gamma ^{\prime }\rangle -2\langle \xi
,\gamma ^{\prime }\rangle \langle J\gamma ^{\prime
},V\rangle=const.
$$
\endproclaim
\demo{Proof} Taking the inner product of (3.8) with $\gamma
^{\prime }$, we get
$$
\frac{d^2}{dt^2}\langle V,\gamma ^{\prime }\rangle =2\langle \xi
,\gamma ^{\prime }\rangle \frac d{dt}\langle J\gamma ^{\prime
},V\rangle.\tag{3.9}
$$
Note that $\langle \xi ,\gamma ^{\prime }\rangle $ is constant
along $\gamma $. It follows from (3.9) that
$$
\frac d{dt}\langle V,\gamma ^{\prime }\rangle =2\langle \xi
,\gamma ^{\prime }\rangle \langle J\gamma ^{\prime },V\rangle
+\alpha
$$
for some constant $\alpha$. \qed
\enddemo

\proclaim{Theorem 3.2} Let $M$ be a Sasakian manifold and let
$\gamma $ be a $\nabla $-geodesic parameterized by arc length.
Then every Jacobi field $V$ along $\gamma $ can be uniquely
decomposed in the following form:
$$
V=a\gamma ^{\prime }+bt\gamma ^{\prime }+W
$$
where $a,b\in R$ and $W$ is a Jacobi field along $\gamma $ such
that
$$
\langle W,\gamma ^{\prime }\rangle =2\langle \xi ,\gamma ^{\prime
}\rangle \int_0^t\langle V,J\gamma ^{\prime }\rangle .
$$
In particular, if either i) $\gamma $ is horizontal, or ii)
$V_{\gamma (t)}\bot J\gamma ^{\prime }(t)$ for every $t$, then $W$
is perpendicular to $\gamma $.
\endproclaim
\demo{Proof} Set
$$
a=\langle V,\gamma ^{\prime }\rangle _{\gamma (0)},\quad b=\langle
V^{\prime },\gamma ^{\prime }\rangle _{\gamma (0)}-2\langle \xi
,\gamma ^{\prime }\rangle \langle J\gamma ^{\prime },V\rangle
_{\gamma (0)}\tag{3.12}
$$
and
$$
W=V-a\gamma ^{\prime }-bt\gamma ^{\prime }.\tag{3.13}
$$
Since $(a+bt)\gamma ^{\prime }$ is a Jacobi field, $W$ is a Jacobi
field too. Then Lemma 3.1 implies that
$$
\frac d{dt}\langle W,\gamma ^{\prime }\rangle =2\langle \xi
,\gamma ^{\prime }\rangle \langle J\gamma ^{\prime },W\rangle
+\beta\tag{3.14}
$$
for some $\beta \in R$. In particular, by taking $t=0$, then
(3.12), (3.13) and (3.14) lead to
$$
\aligned \beta &=\langle W^{\prime },\gamma ^{\prime }\rangle
_{\gamma (0)}-2\langle \xi ,\gamma ^{\prime }\rangle \langle
J\gamma ^{\prime },W\rangle _{\gamma
(0)} \\
&=\langle V^{\prime },\gamma ^{\prime }\rangle _{\gamma
(0)}-b-2\langle \xi ,\gamma ^{\prime }\rangle \langle J\gamma
^{\prime },V\rangle _{\gamma (0)}\\
&=0.\endaligned\tag{3.15}
$$
Note that $\langle W,\gamma ^{\prime }\rangle _{\gamma
(0)}=\langle V,\gamma ^{\prime }\rangle _{\gamma (0)}-a=0$. Then
we integrate (3.14) from $0$ to $t$ and employ (3.15) to find
$$
\langle W,\gamma ^{\prime }\rangle _{\gamma (t)}=2\langle \xi
,\gamma ^{\prime }\rangle \int_0^t\langle J\gamma ^{\prime
},W\rangle dt.
$$
This completes the proof.\qed
\enddemo
\remark{Remark 3.1} \newline (i) If $\gamma $ is a vertical
geodesic, then $J\gamma ^{\prime }=0$. As a result of Theorem 3.2,
any Jacobi field $V$ along the vertical geodesic can be written
uniquely as
$$
V=a\gamma ^{\prime }+bt\gamma ^{\prime }+W
$$
with $\langle W,\gamma ^{\prime }\rangle \equiv 0$. On the other
hand, for any $\nabla $-geodesic in a Sasakian manifold, we have
already shown that $\xi _{\gamma (t)}$ is a Jacobi field which
obviously satisfies $\xi \perp J\gamma ^{\prime }$. Since $\langle
\xi ,\gamma ^{\prime }\rangle $ is constant, we may write $\xi
_{\gamma (t)}=\langle \xi ,\gamma ^{\prime }\rangle \gamma
^{\prime }+W$, where $W=\xi -\langle \xi ,\gamma ^{\prime }\rangle
\gamma ^{\prime }$ is clearly a Jacobi field with $\langle
W,\gamma ^{\prime }\rangle \equiv 0$. \newline (ii) In [BD], the
authors established a similar decomposition for Jacobi fields
along horizontal $\nabla $-geodesics in pseudo-Hermitian
manifolds. Here we give the decomposition for Jacobi fields along
general $\nabla$ -geodesics in Sasakian manifolds.
\endremark

\example{Example 3.1} Let us consider the Jacobi fields along a
$\nabla $-geodesic $\gamma $ in the Heisenberg group $H_m$. Since
$H_m$ is a Sasakian manifold with zero curvature (cf. [DT]), the
Jacobi equation becomes
$$
V^{\prime \prime }=2\langle J\gamma ^{\prime },V^{\prime }\rangle
\xi.\tag{3.16}
$$
First, we assume that $\gamma ^{\prime }$ is not vertical. Set
$\gamma _H^{\prime }=\gamma ^{\prime }-\langle \gamma ^{\prime
},\xi \rangle \xi $. Since $\langle \gamma ^{\prime },\xi \rangle
=$const., we find that $\gamma _H^{\prime }$ and $J\gamma
_H^{\prime }$ are parallel along $\gamma $ and
$$
\Vert\gamma _H^{\prime }\Vert=\Vert J\gamma ^{\prime
}\Vert=\sqrt{1-\langle \gamma ^{\prime },\xi \rangle ^2}\not =0.
$$
Choose a basis $\{v_1,...,v_{2m}\}$ in $H_{\gamma (0)}(H_m)$ such
that
$$
\{v_1=J\gamma _H^{\prime }(0)/\Vert J\gamma _H^{\prime
}(0)\Vert,v_2=\gamma _H^{\prime }(0)/\Vert\gamma _H^{\prime
}(0)\Vert,v_3,...,v_{2m}\}
$$
is an orthonormal basis of $HM_{\gamma (0)}$. Let
$\{E_A(t)\}_{A=0}^{2m}$ be parallel vector fields along $\gamma
(t)$ such that $E_0(t)=\xi _{\gamma (t)} $ and $E_i(0)=v_i$
($i=1,...,2m$). Suppose $V(t)$ is a Jacobi field along $\gamma $.
Then we write
$$
V(t)=\sum_{A=0}^{2m}f_A(t)E_A(t)
$$
and substitute $V$ into (3.16) to find
$$
f_0^{\prime \prime }=2\sqrt{1-\langle \gamma ^{\prime },\xi
\rangle ^2} f_1^{\prime },\quad f_i^{\prime \prime }=0\quad (1\leq
i\leq 2m).\tag{3.17}
$$
Consequently
$$
f_0(t)=a_1\sqrt{1-\langle \gamma ^{\prime },\xi \rangle ^2}
t^2+a_0t+b_0,\quad f_i(t)=a_it+b_i\quad (1\leq i\leq
2m),\tag{3.18}
$$
where $a_A,b_A$ are constants\ ($0\leq A\leq 2m$). Hence we deduce
$$
V=(a_1\sqrt{1-\langle \gamma ^{\prime },\xi \rangle
^2}t^2+a_0t+b_0)\xi _{\gamma
(t)}+\sum_{i=1}^{2m}(a_it+b_i)E_i(t).\tag{3.19}
$$
Next we assume that $\gamma ^{\prime }(t)=\xi _{\gamma (t)}$, and
choose an orthonormal basis $\{v_1,...,v_{2m}\}$ in $H_{\gamma
(0)}M$. Let $\{E_A(t)\}_{A=0}^{2m}$ be parallel vector fields with
$E_0=\xi _{\gamma (t)}$ and $E_i(0)=v_i$ ($i=1,...,2m$). Note that
$J\gamma ^{\prime }=0$ in this case. Similarly we may deduce from
(4.16) the following general solution of the Jacobi equation
$$
V(t)=\sum_{A=0}^{2m}(a_At+b_A)E_i(t).\tag{3.20}
$$
\endexample

\heading{\bf 4. Cartan-Hadamard Type Theorem}
\endheading

We say that a point $q$ is conjugate to $p$ if $q$ is a singular
value of $\exp _p^\nabla :T_pM\rightarrow M$. The conjugacy is
said to be along $ \gamma _v$ if $d\exp _p^\nabla $ is singular at
$v$. The following result is known for any linear connection.

\proclaim{Proposition 4.1} (cf. [KN2]) Let $\gamma
:[0,1]\rightarrow M$ be a $\nabla $-geodesic with $ \gamma (0)=p$,
$\gamma (1)=q$. Then $q$ is conjugate to $p$ along a $\nabla $
-geodesic $\gamma $ if and only if there exists a non-zero Jacobi
field $V$ along $\gamma $ such that $V(0)=V(1)=0$. Hence $q$ is
conjugate to $p$ if and only if $p$ is conjugate to $q$.
\endproclaim
It follows immediately from Proposition 4.1 that
\proclaim{Corollary 4.2} Let $\gamma :[0,l]\rightarrow M$ be a
$\nabla $-geodesic. If $\gamma (0)$ and $\gamma (l)$ are not
conjugate, then a Jacobi field $V$ along $\gamma $ is determined
by its values at $\gamma (0)$ and $\gamma (l)$.
\endproclaim

Let us come back to investigate the conjugate points of
pseudo-Hermitian manifolds. First we consider the special case
that the geodesic is vertical.

\proclaim{Proposition 4.3} Let $(M^{2m+1},H(M),J,\theta )$ be a
Sasakian manifold. If $\gamma :[0,l]\rightarrow M$ is a vertical
geodesic, then $\gamma |_{(0,l]}$ contains no point conjugate to
$\gamma (0)$.
\endproclaim
\demo{Proof} Suppose $\gamma $ is a vertical geodesic and $V$ is a
Jacobi field along $\gamma $ with
$$
V_{\gamma (0)}=V_{\gamma (l)}=0.\tag{4.1}
$$
Since $M$ is Sasakian and $\gamma ^{\prime }(t)=\xi _{\gamma
(t)}$, we have $R(\gamma ^{\prime },V)\gamma ^{\prime }=0$ and
$J\gamma ^{\prime }=0$. Therefore (3.8) becomes $V^{\prime \prime
}=0$. By Corollary 4.2 and (4.1), we conclude that $V\equiv 0$.
\qed
\enddemo

Sometimes it is convenient to consider the decomposition
$$
TM_{\gamma (t)}=span\{\xi \}_{\gamma (t)}\oplus HM_{\gamma
(t)}\tag{4.2}
$$
along the geodesic $\gamma $. Accordingly, we write a vector field
$V$ along $\gamma $ as
$$
V=V_\xi +V_H\tag{4.3}
$$
where $[\cdot ]_H$ and $[\cdot ]_\xi $ denote the horizontal and
vertical components of the vector respectively. Since the
Tanaka-Webster connection $\nabla $ preserves the above
decomposition, the Jacobi equation (3.8) may be decomposed as
$$
\aligned
&\nabla _{\gamma ^{\prime }}\nabla _{\gamma ^{\prime
}}V_H=R(\gamma ^{\prime
},V_H)\gamma ^{\prime } \\
&\nabla _{\gamma ^{\prime }}\nabla _{\gamma ^{\prime }}V_\xi
=2\langle J\gamma ^{\prime },\nabla _{\gamma ^{\prime
}}V_H\rangle\xi
\endaligned\tag{4.4}
$$
by employing the curvature property (1.18) of Sasakian manifolds.
To solve the first equation in (4.4), we may assume the initial
conditions $V_H(0)$ and $V_H^{\prime }(0)$. Whenever $V_H$ is
known, the second equation in (4.4) can be solved by assuming
$V_\xi (0)$ and $V_\xi ^{\prime }(0)$.

The first equation of (4.4) yields
$$
\aligned &\langle V_H,\gamma ^{\prime }\rangle ^{\prime \prime
}=\langle V_H^{\prime \prime },\gamma ^{\prime }\rangle \\
&=\langle R(\gamma ^{\prime },V_H)\gamma ^{\prime },\gamma
^{\prime }\rangle\\
&=0.\endaligned
$$
Therefore the horizontal component $V_H$ may be written uniquely
as
$$
V_H=(at+b)\gamma ^{\prime }+V_H^{\perp }\tag{4.5}
$$
where $\langle V_H^{\perp },\gamma ^{\prime }\rangle \equiv 0$.
Consequently we see that the possible nonlinear part of $\langle
V,\gamma ^{\prime }\rangle $ comes from $\langle V_\xi ,\gamma
^{\prime }\rangle $ which is determined by the integral of
$\langle J\gamma ^{\prime },V_H\rangle $ on $[0,t]$.

The discussion about the expansion $\Vert d\exp ^\nabla
(tw)\Vert^2$ in \S 3 provides some clue to the following result.

\proclaim{Theorem 4.4} Let $(M^{2m+1},H(M),J,\theta )$ be a
Sasakian manifold with $K^H\leq 0$. Let $\gamma :[0,\beta
]\rightarrow M$ be a $\nabla $-geodesic. Then $\gamma ((0,\beta
])$ contains no point conjugate to $\gamma (0)$ along $\gamma $.
Therefore, if $(M,\nabla )$ is complete with $K^H\leq 0$, then $M$
has no conjugate points along any $\nabla $-geodesic.
\endproclaim
\demo{Proof} Suppose $\gamma :[0,l]\rightarrow M$ is a $\nabla
$-geodesic and $V$ is a Jacobi field along $\gamma $ with
$V(0)=V(l)=0$. So $V_H(0)=V_H(l)=0$ and $V_\xi (0)=V_\xi (l)=0$.
From (4.4), we have
$$
\langle \nabla _{\gamma ^{\prime }}\nabla _{\gamma ^{\prime
}}V_H,V_H\rangle =\langle R(\gamma ^{\prime },V_H)\gamma ^{\prime
},V_H\rangle \geq 0.\tag{4.6}
$$
Consequently
$$
\aligned \frac d{dt}\langle \nabla _{\gamma ^{\prime
}}V_H,V_H\rangle &=\langle \nabla _{\gamma ^{\prime }}\nabla
_{\gamma ^{\prime }}V_H,V_H\rangle +\langle
\nabla _{\gamma ^{\prime }}V_H,\nabla _{\gamma ^{\prime }}V_H\rangle \\
&\geq 0,\endaligned\tag{4.7}
$$
that is, $\langle \nabla _{\gamma ^{\prime }}V_H,V_H\rangle (t)$
is non-decreasing. It follow from $V_H(0)=V_H(\beta )=0$ that
$\langle \nabla _{\gamma ^{\prime }}V_H,V_H\rangle \equiv 0$,
which yields
$$
\frac d{dt}\langle V_H,V_H\rangle =2\langle V_H^{\prime
},V_H\rangle =0.\tag{4.8}
$$
From (4.8), we find that $|V_H|=$const., and thus $V_H\equiv 0$.
From the second equation in (4.4), we get
$$
\nabla _{\gamma ^{\prime }}\nabla _{\gamma ^{\prime }}V_\xi =0,
$$
which implies $V_\xi =(at+b)\xi $ for some $a,b\in R$. In view of
$V_\xi (0)=V_\xi (\beta )=0$, we find $V_\xi \equiv 0$. Therefore
we conclude that $V\equiv 0$. \qed
\enddemo
\remark{Remark 4.1} In [BD], the authors proved that if $M$ has
nonpositive horizontal sectional curvature, then there is no
horizontally conjugate point along any horizontal geodesic in a
pseudo-Hermitian manifold.
\endremark

As a consequence of Theorem 4.4, we shall now give a
Cartan-Hadamard type result for Sasakian manifolds with
non-positive horizontal sectional curvature. Let us first recall
the following notion in [DT].

\definition{Definition 4.1} Let $(N,\widetilde{\theta })$ and $(M,\theta )$ be
two pseudo-Hermitian manifolds. We say that a CR map
$f:N\rightarrow M$ is an isopseudo-Hermitian map if $f^{*}\theta=
\widetilde{\theta }$.
\enddefinition

\proclaim{Lemma 4.5} Let $f:(N,\widetilde{\theta })\rightarrow
(M,\theta )$ be an isopseudo-Hermitian map. If $\dim N=\dim
M=2m+1$, then $f:(N,g_{\widetilde{\theta }})\rightarrow
(M,g_\theta )$ is a local isometry with $df(\widetilde{\xi })=\xi
$.
\endproclaim
\demo{Proof} The assumption $f^{*}\theta= \widetilde{\theta }$
implies that $f^{*}d\theta= d\widetilde{\theta }$, and thus
$f^{*}[\theta \wedge (d\theta )^m]=\widetilde{\theta }\wedge
(d\widetilde{\theta })^m$. This yields that $f$ is a local
diffeomorphism. Furthermore, we get
$$
f^{*}G_\theta=\widetilde{G}_{\widetilde{\theta }},\tag{4.9}
$$
since $f$ is a CR map. At any point $p\in N$, we have
$$
\theta_q(df(\widetilde{\xi}_p))=\widetilde{\theta
}_p(\widetilde{\xi}_p)=1\tag{4.10}
$$
where $q=f(p)\in M$. For any vector $Y\in T_qM$, there exists a
vector $X\in T_pN$ such that $df(X)=Y$ , since $df:T_pN\rightarrow
T_qM$ is a linear isomorphism. Therefore
$$
\aligned
d\theta(df(\widetilde{\xi}_p),Y)&=d\theta(df(\widetilde{\xi}_p),df(X)) \\
&=(df^{*}\theta )(\widetilde{\xi}_p,X) \\
&=d\widetilde{\theta }(\widetilde{\xi}_p,X) \\
&=0.
\endaligned\tag{4.11}
$$
Combining (4.10) and (4.11), we find that $df(\widetilde{\xi}
_p)=\xi_q$. Consequently
$$
f^{*}[\theta \otimes \theta +G_\theta]= \widetilde{\theta }\otimes
\widetilde{\theta }+\widetilde{G}_{ \widetilde{\theta }},
$$
that is, $f^{*}g_\theta =g_{\widetilde{\theta }}$.\qed
\enddemo

Next we need the following lemma, whose proof is a slight
modification of that for Lemma 1.32 in [CE].

\proclaim{Lemma 4.6} Let $\varphi :(N,\widetilde{\theta
})\rightarrow (M,\theta )$ be an isopseudo-Hermitian map between
two pseudo-Hermitian manifolds of same dimension. If
$(N,\widetilde{\nabla })$ is complete, then $\varphi $ is a
covering map.
\endproclaim

\demo{Proof} By Lemma 4.5, we know that $\varphi $ is a local
isometry, and it preserves the CR structures. Clearly $\varphi $
maps a $\widetilde{\nabla }$-geodesic to a $\nabla $-geodesic. Fix
$p\in M$ and let $\{p_\alpha \}=\varphi ^{-1}(p)$. Let $D(r)$ and
$D_\alpha (r)$ be the balls about zero of radius $r$ in $T_pM$,
$T_{p_\alpha }N$, respectively. Set $B^\nabla (p;r)=\exp _p^\nabla
(D(r))$ and $B^{\widetilde{\nabla }}(p_\alpha ;r)=\exp _{p_\alpha
}^{\widetilde{\nabla }}(D_\alpha (r))$. Assume $r$ is small enough
so that $B^\nabla (p;r)$ is contained in a normal coordinate
neighborhood around $p$. Write $U=B^\nabla (p;r)$ and $U_\alpha
=B^{\widetilde{\nabla }}(p_\alpha ;r)$ for simplicity. We will
show that $\varphi ^{-1}(U)$ is the disjoint union $\cup _\alpha
U_\alpha $ and that $\varphi :U_\alpha \rightarrow U$ is a
diffeomorphism for each $\alpha $.

Since $\varphi $ preserves locally the pseudo-Hermitian
structures, we have $\varphi (\exp _{p_\alpha }^{\widetilde{\nabla
}})=\exp _p^\nabla d\varphi $, where $d\varphi :T_{p_\alpha
}N\rightarrow T_pM$ is the differential map. Since $\exp _p^\nabla
\circ d\varphi :D_\alpha (r)\rightarrow U$ is a diffeomorphism, so
is $\varphi :U_\alpha \rightarrow U$. Clearly $\cup _\alpha
U_\alpha \subset \varphi ^{-1}(U)$. We shall show the opposite
inclusion. Given $\widetilde{q}\in \varphi ^{-1}(U)$, and set
$q=\varphi (\widetilde{q})$. Let $\gamma $ be the $\nabla
$-geodesic from $q$ to $p$ in the normal coordinate neighborhood
$U$. Let $v=d\varphi ^{-1}(\gamma ^{\prime }(q))\in
T_{\widetilde{q}}N$, and let $\widetilde{\gamma }$ be the $\nabla
$- geodesic from $\widetilde{q}$ in direction $v$. Since
$(N,\widetilde{\nabla })$ is complete, $\widetilde{\gamma }$ may
be extended arbitrary far. Let $t_0\ $be the length of $\gamma $
from $q$ to $p$, and set $\widetilde{p}=\widetilde{\gamma }(t_0)$.
Since $\varphi \circ \widetilde{\gamma }=\gamma $, $\varphi
(\widetilde{p})=p$. Hence $\widetilde{q}\in B^{\widetilde{\nabla
}}(\widetilde{p};r)\subset \cup _\alpha U_\alpha $.

It remains to show that $U_\alpha \cap U_\beta $ is empty if
$\alpha \neq \beta $. Suppose that $U_\alpha \cap U_\beta \neq
\emptyset $. For any $\widetilde{q}\in U_\alpha \cap U_\beta $,
let $\sigma (t)$ (resp. $\tau (t)$) be the geodesic from $p_\alpha
$ to $\widetilde{q}$ in $U_\alpha $ (resp. $p_\mu $ to
$\widetilde{q}$ in $U_\beta $). Then $\varphi (\sigma (t))$ and
$\varphi (\tau (t))$ are two $\nabla $-geodesics from $p$ to
$\varphi (\widetilde{q})$ in $U$. Since $U$ is a normal coordinate
neighborhood, $\varphi (\sigma (t))=\varphi (\tau (t))$, which
implies that $L(\sigma )=L(\tau )$. When $\widetilde{q}$
approaches to $(\partial U_\alpha )\cap U_\beta $, this is
impossible. This completes the proof of the lemma. \qed
\enddemo

\proclaim{Theorem 4.7} Let $(M^{2m+1},H(M),J,\theta )$ be a
Sasakian manifold with $K^H\leq 0$. If $(M,\nabla )$ is complete,
then for any $p\in M$, $\exp _p^\nabla :T_pM\rightarrow M$ is a
covering map. Hence the universal covering space of $M$ is
diffeomorphic to $R^{2m+1}$.
\endproclaim
\demo{Proof} In terms of Theorem 4.4, we find that $\exp _p^\nabla
:T_pM\rightarrow M$ has nonsingular differential. Set
$\widehat{g}=(\exp _p^\nabla )^{*}g_\theta $, $\widehat{\theta
}=(\exp _p^\nabla )^{*}\theta $, $\widehat{H}=\ker \widehat{\theta
}$ and $\widehat{J}=(\exp _p^\nabla )_{*}^{-1}\circ J\circ (\exp
_p^\nabla )_{*}$. Then $(T_pM,\widehat{H},
\widehat{J},\widehat{\theta })$ is a pseudo-Hermitian manifold and
$\exp _p^\nabla :(T_pM,\widehat{\theta })\rightarrow (M,\theta )$
is an isopseudo-Hermitian map.

To simplify notations, we write $N=T_pM$ in what follows. Let
$\widehat{\nabla}$ be the Tanaka-Webster connection of
$(N,\widehat{H}(N),\widehat{\theta},\widehat{J})$. The origin
$o\in T_pM$ is a pole of $(N,\widehat{\nabla })$, since the ray
$\rho (t)=tv$ ($0\leq t<+\infty $) in each direction $v\in T_pM$
is a $\widehat{\nabla }$-geodesic in $N$. First, we claim that
$(N,\widehat{\nabla })$ is complete. Let $\widehat{\sigma
}:[0,r_0)\rightarrow N$ be any $\widehat{\nabla }$-geodesic in $N$
and let $[0,r_0)$ be the largest open interval for which
$\widehat{\sigma }(t)$ exists. Suppose that $r_0<+\infty$. Using
the vector space structure of $N$, we may write $ \widehat{\sigma
}(t)=l(t)v(t)$ with $v(t)\in S^{2m}(1)$, where $l(t)=|
\widehat{\sigma }(t)|_{g_\theta (p)}$ and $S^{2m}(1)$ is the unit
sphere at $ 0$. Set $\gamma _t(s)=\exp _p^\nabla (sv(t))$ ($0\leq
s<+\infty $) and $ \sigma (t)=\exp _p^\nabla (\widehat{\sigma
}(t))$ ($0\leq t<r_0$). Then $\gamma _t$ and $\sigma $ are $\nabla
$-geodesics. Since $(M,\nabla )$ is complete, $\sigma $ may be
extended arbitrarily far. At least, $\sigma :[0,r_0]\rightarrow M$
is well-defined. Let $\{t_i\}$ be a sequence in $[0,r_0]$ such
that $t_i\rightarrow r_0$. By compactness of $S^{2m}(1)$, we may
pass to a sequence of $\{v(t_i)\}$, denoted still by $\{v(t_i)\}$
for simplicity, such that $v(t_i)\rightarrow v$ as $i\rightarrow
+\infty $. Let $ \gamma _v:[0,\infty )\rightarrow M$ be the
$\nabla $-geodesic such that $\gamma _v^{\prime }(0)=v$. Then the
theory of ordinary differential equations (continuous dependence
of solutions on initial data) gives that $\gamma _{t_i}(s)$
converges uniformly to $\gamma _v(s)$ on any closed subinterval of
$[0,\infty )$. Consequently there exists a finite number $l_0$
such $l(t_i)\rightarrow l_0$ and $\gamma _v(l_0)=\sigma (r_0)$.
Set $ \widehat{q}=l_0v\in N$. Let $\widehat{U}$ be a normal
coordinate neighborhood of $\widehat{q}$ as in Theorem 2.1. For
$i$ sufficiently large, $\widehat{\sigma }(t_i)\in \widehat{U}$.
Let $\widehat{\tau }:(-\delta ,\delta )\rightarrow N$ be the
unique $\widehat{\nabla }$-geodesic such that $\widehat{\sigma
}(t_i)\in \widehat{\tau }$ and $\widehat{\tau }(0)=\widehat{q}$.
Then $\widehat{\sigma }\cup \widehat{\tau }$ is a smooth
$\widehat{ \nabla }$-geodesic which extends past $r_0$. This leads
to a contradictions and thus proves the claim. In view of Lemma
4.6, we may conclude that $\exp _p^\nabla :T_pM\rightarrow M$ is a
covering map.\qed
\enddemo

\proclaim{Corollary 4.8} Let $(M^{2m+1},H(M),J,\theta )$ be a
simply connected Sasakian manifold with $K^H\leq 0$. If
$(M,\nabla)$ is complete, then $M$ is diffeomorphic to $R^{2m+1}$.
\endproclaim
\remark{Remark 4.2}
\newline (i) From Theorem 4.7,
we know that any compact Sasakian manifolds with non-positive
horizontal sectional curvature are aspherical. Therefore, by a
standard fact from homotopy theory, their homotopy type is
uniquely determined by their fundamental groups. In [Ch], the
author gave some interesting results about the fundamental groups
of compact Sasakian manifolds.
\newline (ii) We should mention that the above Cartan-Hadamard
type
theorem is invisible from the Levi-Civita connection
$\nabla^{\theta}$ of any Sasakian manifold, because its curvature
tensor with respect to $\nabla^{\theta}$ always satisfies
$R^{\theta}(\xi,X,\xi,X)=1$ for any unit vector $X\in H(M)$.
\endremark

\heading{\bf 5. Basic index lemma and Bonnet-Myers theorem}
\endheading

Let $(M,\theta )$ be a Sasakian manifold and $\gamma
:[0,l]\rightarrow M$ be a $\nabla $-geodesic, parametrized by
arc-length. Given a piecewise differentiable vector field $X$
along $\gamma $, we set
$$
I_0^l(X)=\int_0^l\{\langle \nabla _{\gamma ^{\prime }}X,\nabla
_{\gamma ^{\prime }}X\rangle -\langle R(X,\gamma ^{\prime })\gamma
^{\prime },X\rangle \}dt.\tag{5.1}
$$
The polarization of $I_0^l(X)$ or itself for simplicity will be
called the index form at $\gamma $.

\proclaim{Theorem 5.1} Let $(M,\theta )$ be a Sasakian manifold
and let $\gamma :[0,l]\rightarrow M$ be a $\nabla $-geodesic
parametrized by arc-length, such that $\gamma (0)$ has no
conjugate point along $\gamma $. Let $Y$ be a horizontal Jacobi
field along $\gamma $ such that $Y_{\gamma (0)}=0$ and let $X$ be
a piecewise differentiable vector field along $\gamma $ such that
$X_{\gamma (0)}=0$. If $X_{\gamma (l)}=Y_{\gamma (l)}$ then
$$
I_0^l(X)\geq I_0^l(Y)
$$
and the equality holds if and only if $X=Y$.
\endproclaim
\demo{Proof} Let $J_{\gamma ,0}$ be the space of all Jacobi fields
$Z\in J_\gamma $ such that $Z_{\gamma (0)}=0$. Clearly $\dim
_RJ_{\gamma ,0}=2m+1$. Set $ \widehat{\gamma }_{\gamma
(t)}=t\gamma ^{\prime }$ and $\widehat{\xi } _{\gamma (t)}=t\xi$.
From \S 3, we know that $\widehat{\gamma },\widehat{\xi }\in
J_{\gamma ,0}$.

First we assume that the geodesic $\gamma $ is not vertical. This
implies that $\widehat{\gamma }$ and $\widehat{\xi }$ are linearly
independent in $J_{\gamma ,0}$. Let us complete $\widehat{\gamma
}$ and $\widehat{\xi }$ to a linear basis $\{\widehat{\gamma
},\widehat{\xi },V_2,...,V_{2m}\}$ of $ J_{\gamma ,0}$. Clearly
the vectors $\{\widehat{\gamma }(t),\widehat{\xi }
(t),V_2(t),...,V_{2m}(t)\}$ are linearly independent in $T_{\gamma
(t)}M$ for each $0<t\leq l$, because $\gamma (0)$ has no conjugate
point along $\gamma $.

We set $V_j^H=V_j-\theta (V_j)\xi _{\gamma (t)}$, $2\leq j\leq
2m$, and claim that
$$
\{\widehat{\gamma }(t),\widehat{\xi }_{\gamma (t)},V_{2,\gamma
(t)}^H,...,V_{2m,\gamma (t)}^H\}
$$
are linearly independent for each $t\in (0,l]$ too. To prove this,
suppose that there exist $\{c_A\}_{A=0}^{2m}$ such that
$$
\aligned 0&=c_0\widehat{\gamma }(t)+c_1\widehat{\xi }_{\gamma
(t)}+\sum_{j=2}^{2m}c_iV_{j,\gamma (t)}^H \\
&=c_0\widehat{\gamma }(t)+\{c_1-c_jt^{-1}\theta (V_j)_{\gamma
(t)}\}\widehat{ \xi }_{\gamma (t)}+\sum_{i=2}^{2m}c_jV_j
\endaligned\tag{5.2}
$$
at a fixed $t\in (0,l]$. It follows from (5.2) that
$$
c_0=0,\ c_1-c_jt^{-1}\theta (V_j)_{\gamma (t)}=0,\ c_j=0\quad
(j=2,...,2m),
$$
because $\{\widehat{\gamma }(t),\widehat{\xi }
(t),V_2(t),...,V_{2m}(t)\}$ are linearly independent. Thus $c_A=0$
($A=0,1,...,2m$).

To simplify the notations, let us set $Z_0=\widehat{\gamma
},Z_1=\widehat{\xi }$ and $Z_j=V_j^H$ ($2\leq j\leq 2m$). A simple
observation shows that the vector fields $\{Z_A\}_{A=0}^{2m}$
satisfy the following equation
$$
\nabla _{\gamma ^{\prime }}\nabla _{\gamma ^{\prime }}Z=R(\gamma
^{\prime },Z)\gamma ^{\prime }.\tag{5.3}
$$
Let $\widehat{J}_{\gamma ,0}$ denote the space of all solutions
$Z\ $of (5.3) such that $Z_{\gamma (0)}=0$. Clearly $\dim
_R\widehat{J}_{\gamma ,0}=2m+1$. In fact, $\{Z_A(t)\}_{A=0}^{2m}$
is a linear basis in $\widehat{J}_{\gamma ,0}$ and they are
linearly independent for each $t$. Consequently there are
piecewise differentiable functions $f^A(t)$ such that
$$
X_{\gamma (t)}=\sum_{A=0}^{2m}f^A(t)Z_{A,\gamma (t)}.\tag{5.4}
$$
Now we compute
$$
\aligned |X^{\prime }|^2&=\langle
\sum_A(\frac{df^A}{dt}Z_A+f^AZ_A^{\prime }),\sum_B(
\frac{df^B}{dt}Z_B+f^BZ_B^{\prime })\rangle \\
&=(|\sum_A\frac{df^A}{dt}Z_A|^2+|\sum_Af^AZ_A^{\prime
}|^2)+2\sum_{A,B}\langle \frac{df^A}{dt}Z_A,f^BZ_B^{\prime
}\rangle,\endaligned\tag{5.5}
$$
$$
\aligned -\langle R(X,\gamma ^{\prime })\gamma ^{\prime },X\rangle
&=-\sum_Af^A\langle
R(Z_A,\gamma ^{\prime })\gamma ^{\prime },X\rangle \\
&=\sum_Af^A\langle Z_A^{\prime \prime },X\rangle \\
&=\sum_{A,B}\langle f^AZ_A^{\prime \prime
},f^BZ_B\rangle\endaligned\tag{5.6}
$$
and
$$
\aligned &\sum_{A,B}\langle \frac{df^A}{dt}Z_A,f^BZ_B^{\prime
}\rangle +|\sum_Af^AZ_A^{\prime }|^2+\sum_{A,B}\langle
f^AZ_A^{\prime \prime
},f^BZ_B\rangle \\
&=\frac d{dt}\sum_{A,B}\langle f^AZ_A,f^BZ_B^{\prime }\rangle
-\sum_{A,B}\langle f^AZ_A,\frac{df^B}{dt}Z_B^{\prime }\rangle .
\endaligned\tag{5.7}
$$
Note that
$$
\langle Z_A,Z_B^{\prime }\rangle =\langle Z_B,Z_A^{\prime
}\rangle\tag{5.8}
$$
because $[\langle Z_A,Z_B^{\prime }\rangle -\langle
Z_B,Z_A^{\prime }\rangle ]^{\prime }=0$ and $[\langle
Z_A,Z_B^{\prime }\rangle -\langle Z_B,Z_A^{\prime }\rangle
]_{\gamma (0)}=0$. By employing (5.5), (5.6), (5.7) and (5.8), we
deduce
$$
\aligned &|X^{\prime }|^2-\langle R(X,\gamma ^{\prime })\gamma
^{\prime },X\rangle \\
&=|\sum_A\frac{dfA}{dt}Z_A|^2+\frac d{dt}\sum_{A,B}\langle
f^AZ_A,f^BZ_B^{\prime }\rangle -\sum_{A,B}\langle
f^AZ_A,\frac{df^B}{dt}Z_B^{\prime }\rangle\\
&+\sum_{A,B}\langle \frac{df^A}{dt}Z_A,f^BZ_B^{\prime }\rangle \\
&=|\sum_A\frac{dfA}{dt}Z_A|^2+\frac d{dt}\sum_{A,B}\langle
f^AZ_A,f^BZ_B^{\prime
}\rangle+\sum_{A,B}\frac{df^A}{dt}f^B[\langle Z_A,Z_B^{\prime
}\rangle -\langle
Z_B,Z_A^{\prime }\rangle ] \\
&=|\sum_A\frac{dfA}{dt}Z_A|^2+\frac d{dt}\sum_{A,B}\langle
f^AZ_A,f^BZ_B^{\prime }\rangle .\endaligned\tag{5.9}
$$
Then integrating (5.9) gives
$$
I_0^l(X)=\sum_{A,B}\langle f^AZ_A,f^BZ_B^{\prime }\rangle _{\gamma
(l)}+\int_0^l|\sum_A\frac{df^A}{dt}Z_A|^2\tag{5.10}
$$

Since $Y$ is a horizontal Jacobi field with $Y_{\gamma (0)}=0$, we
know from (4.4) that $Y$ automatically satisfies (5.3). Thus there
are numbers $d^A$ ($A=0,1,...,2m$) such that
$$
Y=\sum_{j=0}^{2m}d^AZ_A.\tag{5.11}
$$
Applying (5.10) to $Y$, we get
$$
I_0^l(Y)=\langle \sum_{A=0}^{2m}a^AZ_A,a^AZ_A^{\prime }\rangle
_{\gamma (l)}.\tag{5.12}
$$
The assumption $X_{\gamma (b)}=Y_{\gamma (b)}$ implies that
$f^A(l)=a^A$ ($0\leq A\leq 2m$). From (5.10) and (5.12), one may
conclude that $I_0^l(X)\geq I_0^l(Y)$ and the equality holds if
and only if $X=Y$.

Now we assume that $\gamma :[0,l]\rightarrow M$ is a vertical
geodesic. In this case, $\widehat{\gamma }=\widehat{\xi }\in
J_{\gamma ,0}$. Let us complete $\widehat{\gamma }$ to a linear
basis $\{\widehat{\gamma } ,V_1,...,V_{2m}\}\in J_{\gamma ,0}$.
Set $Z_0=\widehat{\gamma }$ and $ Z_j=V_j-\theta (V_j)\xi $
($j=1,...,2m$). Since $\gamma (0)$ has no conjugate point along
the vertical geodesic, we may prove similarly that $ \{Z_{0,\gamma
(t)},Z_{1,\gamma (t)},....,Z_{2m,\gamma (t)}\}$ are solutions of
(5.3), and they are linearly independent in $T_{\gamma (t)}M$ for
each $ 0<t\leq l$. Suppose $X$ is a piecewise differentiable
vector field along $\gamma $ such that $X_{\gamma (a)}=0$. Then
$$
X_{\gamma (t)}=\sum_{A=0}^{2m}f^A(t)Z_{A,\gamma (t)}
$$
for some piecewise differentiable functions $f^A$ ($0\leq A\leq
2m$). The remaining argument is similar to the first case. \qed
\enddemo

\remark{Remark 5.1} (i) The above proof actually yields the
following result: Let $Y$ be any solution of (5.3) along $\gamma $
such that $Y_{\gamma (0)}=0$. If $X$ is any piecewise
differentiable vector field along $\gamma $ such that $X_{\gamma
(0)}=0$ and $X_{\gamma (l)}=Y_{\gamma (l)}$, then $I_0^l(X)\geq
I_0^l(Y)$ and the equality holds if and only if $X=Y$. Note that
$Y$ is not necessarily to be horizontal if $Y$ is already a
solution of (6.3); (ii) A similar basic index result was
established in [BD] for a horizontal Jacobi field along a
horizontal $\nabla $-geodesic.
\endremark

Taking $Y=0$ in Theorem 6.1, we have

\proclaim{Corollary 5.2} Let $(M,\theta )$ be a Sasakian manifold
and let $\gamma :[0,l]\rightarrow M$ be a $\nabla $-geodesic,
parametrized by arc length and such that $\gamma (0) $ has no
conjugate point along $\gamma $. If $X$ is a piecewise
differentiable vector field along $\gamma $ such that $X_{\gamma
(0)}=X_{\gamma (l)}=0$, then $I_0^l(X)\geq 0$ and equality holds
if and only if $X=0$.
\endproclaim

As applications of Theorem 5.1, we now prove the Bonnet-Myers type
result for Sasakian manifolds.

\proclaim{Theorem 5.3} Let $(M,\theta )$ be a Sasakian manifold
and let $\gamma :[a,b]\rightarrow M$ be a $\nabla $-geodesic
parametrized by arc length. Assume that the horizontal sectional
curvature satisfies $K^H(\sigma )\geq k_0$ $>0$ for any horizontal
$2$-plane $\sigma \subset T_xM$, $x\in M$. If $b-a\geq \pi /\sqrt{
k_0(1-\langle \gamma ^{\prime },\xi \rangle ^2)}$, then $\gamma
(a)$ has a conjugate point along $\gamma |_{(a,b]}$.
\endproclaim

\demo{Proof} Since there is no conjugate point along a vertical
geodesic, we only need to consider the case that $\gamma ^{\prime
}$ is not parallel to $\xi _{\gamma (t)}$. Let $Y(t)\in
H(M)_{\gamma (t)}$ be a unit vector field along $\gamma $ such
that $\nabla _{\gamma ^{\prime }}Y=0$ and $Y$ is perpendicular to
$\gamma $. Thus $Y\bot \gamma _H^{\prime }$. Set $f(t)=\sin [\pi
(t-a)/(b-a)]$ and $X=f(t)Y$. If $\gamma (a)$ has no conjugate
point in $\gamma |_{(a,b]}$, then Corollary 5.2 implies that
$$
\aligned 0<&I_a^b(X)=\int_a^b\left\{ f^{\prime 2}|Y|^2-f^2\langle
R(Y,\gamma ^{\prime
})\gamma ^{\prime },Y\rangle \right\} dt \\
&=\int_a^b\left\{ f^{\prime 2}|Y|^2-f^2\langle R(Y,\gamma
_H^{\prime })\gamma
_H^{\prime },Y\rangle \right\} dt \\
&\leq \int_a^b\left\{ \frac{\pi ^2}{(b-a)^2}\cos ^2[\pi
(t-a)/(b-a)]-k_0(1-\langle \gamma ^{\prime },\xi \rangle ^2)\sin
^2[\pi
(t-a)/(b-a)]\right\} dt \\
&(s=\pi (t-a)/(b-a) \\
&\leq \frac{b-a}\pi \int_0^\pi \left\{ \frac{\pi ^2}{(b-a)^2}\cos
^2s-k_0(1-\langle \gamma ^{\prime },\xi \rangle ^2)\sin ^2s\right\} ds \\
&=\frac{b-a}2\left\{ \frac{\pi ^2}{(b-a)^2}-k_0(1-\langle \gamma
^{\prime },\xi \rangle ^2)\right\}\endaligned
$$
because $X$ is clearly non-zero. Consequently $b-a<\pi
/\sqrt{k_0(1-\langle \gamma ^{\prime },\xi \rangle ^2)}$. By the
assumption, we conclude that $\gamma (a)$ has a conjugate point in
$\gamma |_{(a,b]}$. \qed

\enddemo

The following result is more general than Theorem 5.3.

\proclaim{Theorem 5.4} Let $(M,\theta )$ be a Sasakian manifold
and let $\gamma :[a,b]\rightarrow M$ be a $\nabla $-geodesic
parametrized by arc length. Assume that the Ricci curvature
satisfies $Ric(X,X)\geq (2m-1)k_0$ $\langle X,X\rangle ,$ $X\in
H(M),$ for some constant $k_0>0$. If $b-a\geq \pi
/\sqrt{k_0(1-\langle \gamma ^{\prime },\xi \rangle ^2)}$, then
$\gamma (a)$ has a conjugate point along $\gamma |_{(a,b]}$.
\endproclaim
\demo{Proof}Suppose that $\gamma :[a,b]\rightarrow M$ is a
geodesic with no conjugate points along $\gamma $ and $\gamma
^{\prime }$ is not parallel to $\xi _{\gamma (t)}$. Let
$\{Y_1,...,Y_{2m-1}\}$ be parallel horizontal vector fields along
$\gamma $ such that $\{\gamma _H^{\prime }/|\gamma _H^{\prime
}|,Y_1,...,Y_{2m-1}\}$ is an orthonormal basis of $H(M)_{\gamma
(t)}$ for every $t$. Set $f(t)=\sin [\pi (t-a)/b-a]$ and
$X_i=f(t)Y_i$ ($i=1,...,2m-1$). Clearly each $X_i$ is a no-zero
vector field with $X_i(\gamma (a))=X_i(\gamma (b))=0$. It follows
from Corollary 5.2 that
$$
\aligned
0<&\sum_{j=1}^{2m-1}I_a^b(X_j)=\sum_{j=1}^{2m-1}\int_a^b\left\{
f^{\prime }(t)^2|Y_j|^2-f^2(t)\langle R(Y_j,\gamma ^{\prime
})\gamma ^{\prime
},Y_j\rangle \right\} \\
&=\int_a^b\left\{ (2m-1)f^{\prime }(t)^2-Ric(\gamma _H^{\prime
},\gamma
_H^{\prime })f^2(t)\right\} dt \\
&\leq \frac{(2m-1)(b-a)}\pi \int_0^\pi \left\{ \frac{\pi
^2}{(b-a)^2}\cos
^2s-k_0(1-\langle \gamma ^{\prime },\xi \rangle ^2)\sin ^2s\right\} ds \\
&=\frac{(2m-1)(b-a)}2\left\{ \frac{\pi ^2}{(b-a)^2}-k_0(1-\langle
\gamma ^{\prime },\xi \rangle ^2)\right\}\endaligned
$$
This implies that $b-a<\pi /\sqrt{k_0(1-\langle \gamma ^{\prime
},\xi \rangle ^2)}$. We conclude that under the assumption
$b-a\geq \pi /\sqrt{k_0(1-\langle \gamma ^{\prime },\xi \rangle
^2)}$, $\gamma (a)$ should has a conjugate point in $\gamma
|_{(a,b]}$.\qed

\remark{Remark 5.2} In view of Proposition 4.3, we see that the
results for conjugate points in Theorems 5.3 and 5.4 are optimal,
since the length $\pi /\sqrt{k_0(1-\langle \gamma ^{\prime },\xi
\rangle ^2)}\rightarrow\infty$ as the constant $\langle \gamma
^{\prime },\xi \rangle\rightarrow 1$.
\endremark

\enddemo

{\bf Acknowledgments}: The authors are grateful to Professors
Shucheng Chang and Sorin Dragomir for their valuable comments.
They also thank Ping Cheng for helpful conversations.

\vskip 0.5 true cm

\vskip 1 true cm
Yuxin Dong

School of Mathematical Science

and

Laboratory of Mathematics for Nonlinear Science

Fudan University,

Shanghai 200433, P.R. China

\vskip 0.2 true cm yxdong\@fudan.edu.cn

\vskip 1 true cm

Wei Zhang

School of Mathematics,

South China University of Technology,

Guangzhou, 510641, P.R. China

\vskip 0.2 cm

sczhangw\@scut.edu.cn

\enddocument